\documentclass{article}
 \usepackage{amsmath}

\newtheorem{thm}{Theorem}[section]

\newtheorem{prop}[thm]{Proposition}

\newtheorem{lem}[thm]{Lemma}

\newtheorem{rem}[thm]{Remark}

\newtheorem{ex}[thm]{Example}

\newcommand{\be}{\begin{equation}}
\newcommand{\ee}{\end{equation}}
\newcommand{\ben}{\begin{enumerate}}
\newcommand{\een}{\end{enumerate}}
\newcommand{\beq}{\begin{eqnarray}}
\newcommand{\eeq}{\end{eqnarray}}
\newcommand{\beqn}{\begin{eqnarray*}}
\newcommand{\eeqn}{\end{eqnarray*}}

\newcommand{\pa}{\partial}

\newcommand{\pxi}{ {\pa \over \pa x^i}}

\newcommand{\pyj}{ {\pa \over \pa y^j}}

\newcommand{\qed}{\hspace*{\fill}Q.E.D.}  

\begin{document}
\title{A Note On Conformal Vector Fields Of  $(\alpha,\beta)$-Spaces }
\author{Guojun Yang\footnote{Supported by the
National Natural Science Foundation of China (11471226) }}
\date{}
\maketitle
\begin{abstract}
 In this paper, we characterize
conformal vector fields of any (regular or singular)
$(\alpha,\beta)$-space with some PDEs. Further, we show some
properties of conformal vector fields of  a class of singular
$(\alpha,\beta)$-spaces satisfying certain geometric conditions.

{\bf Keywords:}  Conformal vector field, $(\alpha,\beta)$-space,
Douglas metric, Landsberg metric

 {\bf MR(2000) subject classification: }
53B40, 53C60
\end{abstract}

\section{Introduction}

Conformal vector fields play an important role in Finsler
geometry.   Some  problems on $(\alpha,\beta)$-metrics can be
solved by constructing a conformal vector field of a Riemann
metric with certain curvature features. For two conformally
related Finsler metrics on a manifold, their conformal vector
fields coincide (\cite{Y1}).

A vector field $V$ on a manifold $M$ has a complete lifted vector
field $V^c$ on $TM$ (see the definition (\ref{XV}) below). Every
conformal vector field $V$ is associated with a scalar function
$c$ called the conformal factor. If $c$ is a constant, $V$ is said
to be homothetic; if $c=0$, $V$ is said to be Killing. As a
special case of conformal vector fields, homothetic vector fields
have some special properties. For example, Huang-Mo obtain the
relation between the flag curvatures of two Finsler metrics $F$
and $\tilde{F}$, where $\tilde{F}$ is defined by $(F,V)$ under
navigation technique for a homothetic vector field $V$ of $F$
(\cite{MH}).

 An  $(\alpha,\beta)$-metric is defined by
 $$F=\alpha \phi(s),\ \ s=\beta/\alpha,$$
where $\alpha=\sqrt{a_{ij}y^iy^j}$ is a  Riemann metric,
$\beta=b_iy^i$ is a 1-form and $\phi(s)$ is a function satisfying
certain conditions. If taking $\phi(s)=1+s$, we get
$F=\alpha+\beta$, which is called a Randers metric. In \cite{SX},
Shen-Xia study conformal vector fields of (regular) Randers spaces
under certain curvature conditions. In \cite{HM1}, Huang-Mo show
that a conformal vector field of a (regular) Randers space of
isotropic S-curvature must be homothetic.

For a non-Riemannian $(\alpha,\beta)$-metric
$F=\alpha\phi(\beta/\alpha)$ with $\phi(s)$ being $C^{\infty}$ on
an open neighborhood of $s=0$ and $\phi(0)\ne 0$, we characterize
in \cite{Y1} its conformal vector fields by a  systems of PDEs
(cf. \cite{K} \cite{SX}). Further we prove that any conformal
vector field on such an $(\alpha,\beta)$-space satisfying certain
curvature conditions   must be homothetic, and we also give
examples to indicate that a conformal vector field on a
projectively flat Randers space is not necessarily homothetic
(\cite{Y1} \cite{Y2}).

A Finsler metric $F>0$ on a manifold $M$ is said to be regular if
$F$ is positively definite on the whole slit tangent bundle
$TM-0$. Otherwise, $F (>0)$ is said to be singular. Singular
Finsler metrics have a lot of applications in the real world
(\cite{AIM} \cite{AHM}). Z. Shen also introduces singular Finsler
metrics in \cite{Shen3}. For an $(\alpha,\beta)$-metric $F$
discussed in this paper, we do not assume that $F$ is regular, but
 $\alpha$ is supposed to be regular (sometimes $\alpha$ can even be singular). If
$\phi(0)$ is not defined or $\phi(0)=0$,
 then the $(\alpha,\beta)$-metric $F= \alpha \phi
(\beta/\alpha)$ is singular.  Two $(\alpha,\beta)$-metrics $F$ and
$\widetilde{F}$ are said to be of the same metric type if they can
be written as
 $$F=\alpha\phi(\beta/\alpha),\
 \widetilde{F}=\widetilde{\alpha}\phi(\widetilde{\beta}/\widetilde{\alpha}):\
 \ \
 \widetilde{\alpha}=\sqrt{k_1\alpha^2+k_2 \beta^2},\
\widetilde{\beta}=k_3\beta,$$
 where $k_1(>0),k_2,k_3(\ne 0)$ are constant. We are mainly concerned about
 two singular metric types: $m$-Kropina type $F=\beta^m\alpha^{1-m}$ ($m\ne
 0,1$) and the  type
 $F=\beta e^{\epsilon\alpha^2/\beta^2}$ ($\epsilon=\pm
 1$).  We put a set
  \be\label{y1}
 \Theta:= \big\{F=\beta^m\alpha^{1-m} \ (m\ne
 0,1),\ \  F=\beta e^{\pm\alpha^2/\beta^2}\big\}.
  \ee
  The case $m=-1$ is called a Kropina metric and it is first
  studied by Kropina in \cite{Kr}.

In this paper, we will characterize conformal vector fields of any
(regular or singular) $(\alpha,\beta)$-space, and investigate some
properties of  conformal vector fields of the singular
$(\alpha,\beta)$-spaces (\ref{y1}) satisfying certain conditions.

\begin{thm}\label{th1}
Let $F= \alpha \phi (\beta/\alpha)$ be a non-Riemannian
$(\alpha,\beta)$-metric being not of the metric type in $\Theta$
(see (\ref{y1})). Then $V$ is a conformal vector field of $F$ with
the conformal factor $c$ if and only if $V$ satisfies
 \be\label{y2}
  V^c(\alpha^2)=4c\alpha^2,\ \ \ \ V^c(\beta)=2c\beta.
  \ee
\end{thm}

For an $m$-Kropina metric $F=\beta^m\alpha^{1-m} \ (m\ne
 0,1)$, we can always put $||\beta||_{\alpha}=1$ without loss of generality (see
  Lemma \ref{lem42} below) (cf. \cite{SY} \cite{Y3} \cite{Y4}).  Theorem
 \ref{th1} generalizes the corresponding result in \cite{Y1} for regular $(\alpha,\beta)$-metrics. In Theorem
 \ref{th1}, if $F$ is a metric listed in (\ref{y1}), we have the
 following characterization result for conformal vector fields.

\begin{thm}\label{th2}
 Let $F$ be an $(\alpha,\beta)$-metric defined by (\ref{y1}) on a manifold $M$ and
 $V$ be a vector field on $M$. If $F=\beta^m\alpha^{1-m}$ is an
  $m$-Kropina metric with $||\beta||_{\alpha}=1$,
  then $V$ is a conformal vector
field of $F$ with the conformal factor $c$ iff. $V$ satisfies
(\ref{y2}), namely,
 \be\label{y3}
  V^c(\alpha^2)=4c\alpha^2,\ \ \ \ V^c(\beta)=2c\beta.
  \ee
  If $F$ is given by $F=\beta e^{\pm\alpha^2/\beta^2}$, then $V$ is a conformal vector
field of $F$ with the conformal factor $c$ iff. $V$ satisfies
 \be\label{y4}
 V^c(\alpha^2)=2\tau \alpha^2\pm(2c-\tau)\beta^2,\ \ \ \ V^c(\beta)=\tau
 \beta,
 \ee
 where $\tau$ is a scalar function on $M$.
\end{thm}

Using some basic results in \cite{Y2} and Theorem \ref{th2}, we
 further in Section \ref{sec5} study
conformal vector fields of the  $(\alpha,\beta)$-metrics listed in
(\ref{y1}) under some  conditions. For an $m$-Kropina metric, we
get Theorems \ref{th51} and \ref{th52} below. For a Kropina metric
of Douglas type, we construct Example \ref{ex1} below to show the
existence of non-homothetic conformal vector fields. For the
second metric in (\ref{y1}), we obtain Theorem \ref{th53} below.

\section{Preliminaries}

Let $F$ be a Finsler metric on  a manifold $M$, and $V$ be a
vector field on $M$. Let $\varphi_t$ be the flow generated by $V$.
Define $\widetilde{\varphi}_t: TM \mapsto TM$ by
  $\widetilde{\varphi}_t(x,y)=(\varphi_t(x),\varphi_{t*}(y))$.
  $V$ is said to be conformal if (cf. \cite{HM1})
 \be\label{y5}
 \widetilde{\varphi}_t^*F=e^{2\sigma_t}F,
 \ee
where $\sigma_t$ is a function on $M$ for every $t$.
 Differentiating (\ref{y5}) by $t$ at $t=0$,  we obtain
  $$
 V^c(F)=2cF,
  $$
  where we define
   \be\label{XV}
   V^c:=V^i\pxi+y^i\frac{\pa
   V^j}{\pa x^i}\pyj,\ \ \ \ \  c:=\frac{d}{dt}|_{t=0}\sigma_t.
   \ee
 In (\ref{XV}), the function $c$ is called the conformal factor,
 and $V^c$ is called the complete lift of $V$.

\begin{lem}\label{lem20}
A vector
   field $V$ on a Finsler manifold $(M,F)$ is conformal with the conformal factor $c$  if and
   only if
    $$
  V^c(F^2)=4cF^2  \ (\Longleftrightarrow
 V^c(F)=2cF),\ \ \  or \  V_{0|0}=2cF^2.
   $$
   where $_|$ is the $h$-covariant derivative  of Cratan (Berwald, or Chern) connection.
\end{lem}

\begin{lem}\label{lem21}
Let $\beta=b_i(x)y^i$ be a 1-form, and $V$ be a vector field on a
Riemann manifold $(M,\alpha)$ with $\alpha=\sqrt{a_{ij}y^iy^j}$.
Then we have
 \be\label{y7}
 V^c(\alpha^2)=2V_{0|0},\ \ \ \ \ V^c(\beta)=(V^j\frac{\pa b_i}{\pa
 x^j}+b_j\frac{\pa V^j}{\pa x^i})y^i=(V^jb_{i|j}+b^jV_{j|i})y^i,
 \ee
 where $V_i:=a_{ij}V^j$ and $b^i:=a^{ij}b_j$, and the covariant derivative is taken with respect to the
 Levi-Civita connection of $\alpha$.
\end{lem}

In this paper, for a Riemannian metric $\alpha
=\sqrt{a_{ij}y^iy^j}$ and a $1$-form $\beta = b_i y^i $,  let
 $$r_{ij}:=\frac{1}{2}(b_{i|j}+b_{j|i}),\ \ s_{ij}:=\frac{1}{2}(b_{i|j}-b_{j|i}), \ \ s_j:=b^is_{ij},
 \ \   b:=||\beta||_{\alpha},$$
 where we define $b^i:=a^{ij}b_j$, $(a^{ij})$ is the inverse of
 $(a_{ij})$, and $\nabla \beta = b_{i|j} y^i dx^j$  denotes the covariant
derivatives of $\beta$ with respect to $\alpha$.

\section{Proof of Theorem \ref{th1}}

By Lemma \ref{lem20} ,  $V$ is a conformal vector field of $F$
with the conformal factor $c$ iff. $V^c(F^2)=4cF^2$. A direct
computation shows that
 \beqn
 V^c(F^2)&=&\phi^2V^c(\alpha^2)+2\alpha^2\phi\phi'\frac{\alpha
 V^c(\beta)-\beta
 V^c(\alpha)}{\alpha^2}\\
 &=&\phi(\phi-s\phi')V^c(\alpha^2)+2\alpha\phi\phi'V^c(\beta).
 \eeqn
Now plugging (\ref{y7}) into the above equation, we see that
$V^c(F^2)=4cF^2$ is written as
 \be\label{y8}
 V_{0;0}+\alpha
 Q(V^ib_{j;i}+b^iV_{i;j})y^j=\frac{2c\phi}{\phi-s\phi'}\alpha^2,\
 \ \ \  (Q:=\frac{\phi'}{\phi-s\phi'}),
  \ee
where the covariant derivative is taken with respect to
  $\alpha$.

In order to simplify (\ref{y8}), we choose a special coordinate
 system $(s,y^a)$ at  a fixed point on a manifold as usually used.
Fix an arbitrary point $x\in M$ and take  an orthogonal basis
  $\{e_i\}$ at $x$ such that
   $$\alpha=\sqrt{\sum_{i=1}^n(y^i)^2},\ \ \beta=by^1.$$
It follows from $\beta=s\alpha$ that
 $$y^1=\frac{s}{\sqrt{b^2-s^2}}\bar{\alpha},\ \ \ \ \Big(\bar{\alpha}:=\sqrt{\sum_{a=2}^n(y^a)^2}\Big).$$
Then if we  change coordinates $(y^i)$ to $(s, y^a)$, we get
  $$\alpha=\frac{b}{\sqrt{b^2-s^2}}\bar{\alpha},\ \
  \beta=\frac{bs}{\sqrt{b^2-s^2}}\bar{\alpha}.$$
 Let
 $$\bar{V}_{0;0}:=V_{a;b}y^ay^b, \ \ \bar{V}_{1;0}:=V_{1;a}y^a,\ \ \bar{V}_{0;1}:=V_{a;1}y^a,
 \ \ \bar{b}_{0;i}:=b_{a;i}y^a $$
Note that under the coordinate $(s,y^a)$, we have
$b_1=b,\bar{b}_0=0$, but generally $\bar{b}_{0;i}\ne 0$.

Under the coordinate $(s,y^a)$, (\ref{y8}) is equivalent to
  \beq
 0&=&b(b\bar{V}_{1;0}+V^i\bar{b}_{0;i})Q+
 (\bar{V}_{1;0}+\bar{V}_{0;1})s,\label{y9}\\
 0&=&\big[b(-2bc+V^ib_{1|i}+bV_{1|1})sQ-2b^2c+V_{1|1}s^2\big]\bar{\alpha}^2+(b^2-s^2)\bar{V}_{0|0}.\label{y10}
  \eeq

\begin{lem}\label{lem31}
 If $Q\ne k_1s+k_2/s$ for any constants $k_1$ and $k_2$, then (\ref{y9}) and (\ref{y10}) are
 equivalent to
 \beq
&&b\bar{V}_{1;0}+V^i\bar{b}_{0;i}=0,\ \ \ \
 \bar{V}_{1;0}+\bar{V}_{0;1}=0,\label{y11}\\
&&\bar{V}_{0;0}=2c\bar{\alpha}^2,\ \ \ V^ib_{1|i}+bV_{1|1}=2bc,\ \
\  V_{1|1}=2c.\label{y12}
 \eeq
\end{lem}

{\it Proof :} We only need to show that (\ref{y9}) and (\ref{y10})
imply (\ref{y11}) and (\ref{y12}). If
$b\bar{V}_{1;a}+V^i\bar{b}_{a;i}\ne 0$ (for some $a\ge 2$) at a
point, then by (\ref{y9})  we have $Q=k_1s$, where
 $$
 k_1:=
 -\frac{\bar{V}_{1;a}+\bar{V}_{a;1}}{b(b\bar{V}_{1;a}+V^i\bar{b}_{a;i})},\
 \ (a\ge 2).
 $$
 Since $Q=Q(s)$ and $k_1=k_1(x)$, we see that $k_1$ is a constant. Then it is a contradiction by assumption. So we
have the first equation in (\ref{y11}). Then it is easy to get the
second equation in (\ref{y11}) by (\ref{y9}).

 By (\ref{y10}), we first get $\bar{V}_{0;0}=2\tau\bar{\alpha}^2$
 for a scalar function $\tau=\tau(x)$. Then plugging it into
 (\ref{y10}) we have
  \be\label{y13}
b(-2bc+V^ib_{1|i}+bV_{1|1})sQ+(V_{1|1}-2\tau)s^2+2(\tau-c)b^2=0.
  \ee
  If $-2bc+V^ib_{1|i}+bV_{1|1}\ne 0$ at a point, then by
  (\ref{y13}) we have
  \be\label{y14}
Q=k_1s+\frac{k_2}{s},
  \ee
  where $k_1$ and $k_2$ are defined by
   $$
 k_1:=-\frac{V_{1|1}-2\tau}{b(-2bc+V^ib_{1|i}+bV_{1|1})},\ \
 k_2:=-\frac{2(\tau-c)b^2}{b(-2bc+V^ib_{1|i}+bV_{1|1})},
   $$
  which can be proved to be constant  by a similar analysis as the
  above. By assumption, we get the second equation in (\ref{y12}).
  Now it is easy to get $V_{1|1}=2\tau,\tau=c$ by (\ref{y13}).
  Therefore, all equations in (\ref{y12}) hold. \qed

\

Solving the ODE (\ref{y14}), we have two cases: (i) if $k_2=-1$
($k_1\ne 0$), then
 \be\label{y15}
 \phi(s)= c_1se^{\frac{c_2}{s^2}},\ \  (c_2:=1/(2k_1));
 \ee
 (ii) if $k_2\ne -1$, then
 \be\label{y16}
 \phi(s)=c_1(1+c_2s^2)^{\frac{1-m}{2}}s^m,\ \
 (m:=k_2/(1+k_2),\ \ c_2:=k_1/(1+k_2)),
 \ee
where $c_1$ is a constant. Now if $\phi(s)$ is given by
(\ref{y15}) or (\ref{y16}), then the $(\alpha,\beta)$-metric $F$
is Riemannian (if $m=0$ in (\ref{y16})), or is of the metric type
listed in (\ref{y1}). So by the assumption of Theorem \ref{th1},
we have $Q\ne k_1s+k_2/s$ for any constants $k_1$ and $k_2$. Thus
we obtain (\ref{y11}) and (\ref{y12}) by Lemma \ref{lem31}. Under
arbitrary coordinate system, (\ref{y11}) and (\ref{y12}) are
written as
 $$
V_{i|j}+V_{j|i}=4c a_{ij},\ \ \  V^jb_{i|j}+b^jV_{j|i}=2c b_i,
 $$
where the covariant derivative is taken with respect to $\alpha$.
The above equations are equivalent to (\ref{y2}) by Lemma
\ref{lem21}. This completes the proof of Theorem \ref{th1}. \qed

\section{Proof of Theorem \ref{th2}}\label{sec4}

\subsection{The metric of $m$-Kropina type}

\begin{prop}\label{prop41}
 Let $F=(\alpha^2+k\beta^2)^{(1-m)/2}\beta^m$ be a  metric of
$m$-Kropina type.  Then $V=V^i\pa/\pa x^i$ is a conformal vector
field of $F$ with the conformal factor $c$ if and only if $V$
satisfies
 \be\label{y17}
  V_{i|j}+V_{j|i}=2\tau a_{ij}-\frac{2k(2c-\tau)}{m}b_ib_j,\ \ \ \
  V^jb_{i|j}+b^jV_{j|i}=\big(\tau+\frac{2c-\tau}{m}\big)b_i,
  \ee
where $\tau=\tau(x)$ is a scalar function and the covariant
derivatives are taken with respect to the Levi-Civita connection
of $\alpha$.
\end{prop}

{\it Proof :} $V=V^i\pa/\pa x^i$ is a conformal vector field of
$F$ with the conformal factor $c$ iff. (\ref{y8}) holds, where
$\phi(s)=(1+ks^2)^{(1-m)/2}s^m$. It is easy to show that
(\ref{y8}) is equivalent to
 $$
 \big[2c\beta-m(b^iV_{i|0}+V^i\beta_{|i})\big]\alpha^2+\beta\big[(m-1)V_{0|0}+2kc\beta^2-k\beta(b^iV_{i|0}+V^i\beta_{|i})\big]=0.
 $$
It is easy to see that, for some scalar function $\tau=\tau(x)$,
the above equation  is equivalent to
  \beq
&&(m-1)V_{0|0}+2kc\beta^2-k\beta(b^iV_{i|0}+V^i\beta_{|i})=(m-1)\tau\alpha^2,\label{y18}\\
&& 2c\beta-m(b^iV_{i|0}+V^i\beta_{|i})=-(m-1)\tau\beta,\label{y19}
  \eeq
 Rewriting (\ref{y18}) and (\ref{y19}), we
immediately obtain (\ref{y17}).  \qed

\

In Proposition \ref{prop41}, letting $k=0$, we immediately obtain
the characterization equations for the conformal vector field of
an $m$-Kropina metric. Due to the following known lemma, we can
always put $b=||\beta||_{\alpha}=1$ without loss of generality.

\begin{lem} \label{lem42} (\cite{SY} \cite{Y3} \cite{Y4})
For and $m$-Kropina metric $F=\beta^m\alpha^{1-m}$, we have
$$
F=\beta^m\alpha^{1-m}=\widetilde{\beta}^m\widetilde{\alpha}^{1-m},\
\ \  (\widetilde{\alpha}:=b^m\alpha, \
\widetilde{\beta}:=b^{m-1}\beta).
$$
\end{lem}

{\it Proof of Theorem \ref{th2} :}  For an $m$-Kropina metric
$F=\beta^m\alpha^{1-m}$ with $b=||\beta||_{\alpha}=1$, let $V$ be
a conformal vector field of $F$ with the conformal factor $c$.
Then by (\ref{y17}) we have
 \be\label{y20}
V_{i|j}+V_{j|i}=2\tau a_{ij},\ \ \ \
  V^jb_{i|j}+b^jV_{j|i}=\big(\tau+\frac{2c-\tau}{m}\big)b_i.
 \ee
Using $b=1$ and contracting the second equation of (\ref{y20}) by
$b^i$, we get
 \be\label{y21}
 b^ib^jV_{j|i}=\tau+\frac{2c-\tau}{m}, \ \ (since \ b^ib_{i|j}=0).
 \ee
 By the first equation of (\ref{y20}) we have
 $b^ib^jV_{j|i}=\tau$. Therefore, we get $\tau=2c$ from
 (\ref{y21}). Then by (\ref{y20}) again and $\tau=2c$, we obtain
 (\ref{y3}) from Lemma \ref{lem21}.

 \subsection{The metric $F=\beta e^{\pm\alpha^2/\beta^2}$}

Now we prove  Theorem \ref{th2} when $F$ is the metric $F=\beta
e^{\pm\alpha^2/\beta^2}$. Let $V=V^i\pa/\pa x^i$ be a conformal
vector field of $F$ with the conformal factor $c$. Plugging
$\phi(s)=se^{\pm1 /s^2}$ into (\ref{y8}), we have
 \be\label{y22}
 \pm
 2(b^iV_{i|0}+V^i\beta_{|i})\alpha^2-\beta\big[\beta(b^iV_{i|0}+V^i\beta_{|i})-2c\beta^2\pm
 2V_{0|0}\big]=0.
 \ee
 It is easy to show that (\ref{y22}) is equivalent to
  \be\label{y23}
\beta(b^iV_{i|0}+V^i\beta_{|i})-2c\beta^2\pm
 2V_{0|0}=\pm2\tau\alpha^2,\ \ \
 b^iV_{i|0}+V^i\beta_{|i}=\tau\beta,
  \ee
where $\tau=\tau(x)$ is a scalar function. Now solving
(\ref{y23}), we have
 $$
 V_{i|j}+V_{j|i}=2\tau a_{ij}\pm (2c-\tau)b_ib_j,\ \ \  V^jb_{i|j}+b^jV_{j|i}=\tau b_i,
 $$
 Thus by Lemma \ref{lem21}, the above equations are rewritten as (\ref{y4}).  \qed

\section{Conformal vector fields of the metrics in
(\ref{y1})}\label{sec5}

On the basis of Theorem \ref{th2}, we study some properties of
conformal vector fields of an $m$-Kropina metric or the metric
$F=\beta e^{\pm\alpha^2/\beta^2}$ listed in (\ref{y1}).

We first introduce a result actually  proved in \cite{Y2} as
follows.

\begin{lem} \label{lem51}(cf. \cite{Y2})
 Let $\alpha=\sqrt{a_{ij}y^iy^j}$ be a Riemann metric and
 $\beta=b_iy^i$ be a 1-form and $V=V^i\pa/\pa x^i$  be a vector field on a manifold
 $M$. Suppose $\beta$ is a conformal 1-form of $\alpha$, and
  \be\label{y24}
 V_{i|j}+V_{j|i}=\sigma a_{ij},\ \ \  V^jb_{i|j}+b^jV_{j|i}=\tau b_i,
  \ee
  where $\sigma,\tau$ are scalar functions on $M$, and the covariant derivatives are taken with
respect to the Levi-Civita connection of $\alpha$. Then
$\tau-\sigma$ is a constant.
\end{lem}

\subsection{On $m$-Kropina metrics}

For an $m$-Kropina metric $F=\alpha^{1-m}\beta^m$, we have the
following Theorem.

\begin{thm}\label{th51}
Let $F=\alpha^{1-m}\beta^m$ be an $m$-Kropina metric  and $V$ be a
conformal vector field of $F$. Suppose $F$ is a Landsberg metric
in the dimension $n\ge 3$, or a Douglas metric with $m\ne -1$.
Then $V$ is homothetic.
\end{thm}

{\it Proof :} Let $F$ be a Landsberg metric in the dimension $n\ge
3$, or a Douglas metric with $m\ne -1$. It follows from
\cite{Shen2} (Landsberg case), or \cite{Y3} \cite{Y4} (Douglas
case) that $F$ is characterized by
 \beq
   s_{ij}&=&\frac{b_is_j-b_js_i}{b^2},\label{y26}\\
   r_{ij}&=&2\tau \big\{mb^2a_{ij}-(m+1)b_ib_j\big\}
   -\frac{m+1}{(m-1)b^2}(b_is_j+b_js_i),\label{y27}
 \eeq
where $\tau=\tau(x)$ is a scalar function. By Lemma \ref{lem42},
we put $b=1$ in (\ref{y26}) and (\ref{y27}). Then it is easy to
show that  (\ref{y26}) and (\ref{y27}) are reduced to $b_{i|j}=0$
(cf. \cite{Y3} \cite{Y4}).

Let $V$ be a conformal vector field of $F$ with the conformal
factor $c$. Then we have (\ref{y3}) by Theorem \ref{th2}. So
(\ref{y24}) holds with $\tau=2c,\sigma=4c$, and $c$ is constant by
Lemma \ref{lem51}. Thus $V$ is homothetic. \qed

\begin{rem}
Let $F=\alpha^{1-m}\beta^m$ be an $m$-Kropina metric ($m\ne -1$)
which is of scalar flag curvature in the dimension $n\ge 3$, or a
weak Einstein-metric. Then any conformal vector field  $V$ of $F$
is homothetic. This result follows from \cite{Y5} \cite{SY} which
says that $F$ is flat-parallel, or Ricci-flat parallel, if
$||\beta||_{\alpha}=1$.
\end{rem}

 We prove in \cite{Y3} \cite{Y4} that a
Kropina metric $F=\alpha^2/\beta$ is a Douglas metric if and only
if
 \be\label{y028}
 s_{ij}=\frac{b_is_j-b_js_i}{b^2}.
 \ee
 In Theorem \ref{th51}, if $F=\alpha^2/\beta$ is a Kropina metric
of Douglas type, we will show that its conformal vector fields are
not necessarily homothetic, just like that for Randers metrics of
Douglas type.  In the following, we will give such a family of
examples.

Define $\alpha$ and $\beta$ by
  \beq
 \alpha:&&\hspace{-0.6cm}=\frac{2}{1+\mu|x|^2}|y|, \ \ \ \ \ \
 \beta:=f(c)c_{x^i}y^i,\label{y28}\\
  c:&&\hspace{-0.6cm}=\frac{\tau(1-\mu |x|^2)+ \langle \mu
\gamma+\eta,x\rangle}{1+\mu |x|^2},\label{y29}\\
 V^i:&&\hspace{-0.6cm}=-2\big(\tau+\langle
 \eta,x\rangle\big)x^i+|x|^2\eta^i+q_r^ix^r+\gamma^i,\label{y30}
  \eeq
  where $f$ is a function,  $\mu\ (\ne0)$ and $\tau$ and  $\eta=(\eta^i)$ and
  $\gamma=(\gamma^i)$ are of constant values, and $Q=(q^i_j)$ is a
  constant skew-symmetric matrix, and these  parameters
  satisfy
   \beq
 Q\eta&&\hspace{-0.6cm}=-\mu(4\tau\gamma+Q\gamma), \ \ \ \  \  |\eta|^2=\mu(\mu|\gamma|^2-4\tau^2),
 \label{y31}\\
  f'(c)&&\hspace{-0.6cm}=\frac{2(c-\tau)}{2\tau c-2c^2+\mu|\gamma|^2+\langle
 \eta,\gamma\rangle}f(c).\label{y32}
   \eeq
If we let $F=\alpha+\beta$, then $F$ is projectively flat, and $V$
is a non-homothetic conformal vector field of $F$ (see \cite{Y2}).
To construct a family of Kropina metrics of Douglas type with
non-homothetic conformal vector fields, we put  more conditions on
$\alpha$ and $\beta$ in (\ref{y28}).

Now in (\ref{y28})--(\ref{y32}), put either of the following two
conditions:
 \beq
f^2(c)&&\hspace{-0.6cm}=-\frac{1}{\mu c^2} ,\ \ \ \  \ \ \ \ \
 \langle
\eta,\gamma\rangle=-\mu|\gamma|^2;\label{y34}\\
f^2(c)&&\hspace{-0.6cm}=\frac{2}{\mu(\langle
 \eta,\gamma\rangle+\mu|\gamma|^2-2c^2)},  \ \ \ \   \tau=0.\label{y35}
 \eeq
It is easy to verify that $||\beta||_{\alpha}=1$. Define a Kropina
metric $F=\alpha^2/\beta$. Since $\beta$ is closed, we see from
(\ref{y028}) that $F$ is a Douglas metric.

\begin{ex}\label{ex1}
 Let $\mu,\tau,\eta,\gamma$ and $Q$ satisfy (\ref{y31}), and
 (\ref{y34}) or (\ref{y35}). Define a Kropina
metric $F=\alpha^2/\beta$ by (\ref{y28}), where $c$ is given by
(\ref{y29}) and $f(c)$ is given by (\ref{y34}) or (\ref{y35}). Let
$V$ be vector field given by (\ref{y30}). Then $F$ is a Douglas
metric, and $V$ is a non-homothetic conformal vector field of $F$
if $\eta\ne -\mu\gamma$.
\end{ex}

We have not found a non-homothetic conformal vector field for a
locally projectively flat Kropina metric (cf. \cite{Y2} for
Randers metrics). On the other hand,   for a Kropina metric under
certain curvature conditions, any of its conformal vector fields
is homothetic.

\begin{thm}\label{th52}
 Let $F=\alpha^2/\beta$ be a (weak) Einstein-Kropina metric.
 Then any conformal vector field $V$ of $F$ is homothetic. Further, $V$
 can be locally determined if $F$ is of constant flag curvature.
\end{thm}

We have proved in \cite{SY} that a Kropina metric
$F=\alpha^2/\beta$ is a weak Einstein-metric iff. $F$ is an
Einstein-metric, iff. $r_{00}=0$ and $\alpha$ is an
Einstein-metric (here $||\beta||_{\alpha}=1$). Using $r_{00}=0$ (a
Killing form $\beta$), any conformal vector field in Theorem
\ref{th52} is homothetic by Theorem \ref{th2} and  Lemma
\ref{lem51}. When $F$ is of constant flag curvature in Theorem
\ref{th52}, the local structure of $V$  can be obtained by solving
the equations (\ref{y3}) since $\alpha$ is of constant sectional
curvature and $\beta$ is a Killing form (cf. \cite{Y2}).

\subsection{On the metric $F=\beta e^{\pm\alpha^2/\beta^2}$}

For the metric $F=\beta e^{\pm\alpha^2/\beta^2}$, we will prove
the following theorem.

\begin{thm}\label{th53}
Let $F=\beta e^{\pm\alpha^2/\beta^2}$ be a Landsberg metric in the
dimension $n\ge 3$, or a Douglas metric on the manifold $M$. Then
any conformal vector field of $F$ must be homothetic.
\end{thm}

{\it Proof :} We break the proof of Theorem \ref{th53} into the
following lemmas. First we give a lemma to characterize the metric
$F=\beta e^{\pm\alpha^2/\beta^2}$ which is a Landsberg metric  or
a Douglas metric.

\begin{lem}\label{lem57}
The metric $F=\beta e^{\pm\alpha^2/\beta^2}$ is a Landsberg metric
in the dimension $n\ge 3$, or a Douglas metric if and only if
 \be\label{y39}
 r_{ij}=\sigma \big[(\pm\frac{1}{2}b^2-1)b_ib_j+b^2a_{ij}\big]+\big(\pm 1-\frac{1}{b^2}\big)(b_is_j+b_js_i),\ \ \ \
 s_{ij}=\frac{b_is_j-b_js_i}{b^2},
 \ee
 where $\sigma=\sigma(x)$ is a scalar function. In this case, $F$
 is a Berwald metric.
\end{lem}

{\it Proof :} The Landsberg case has been proved in \cite{Shen2}.
We prove the Douglas case. It is known that an
$(\alpha,\beta)$-metric $F=\alpha\phi(\beta/\alpha)$ is a Douglas
metric iff.
 \beq
 &&\alpha Q (s^i_{\ 0}y^j-s^j_{\ 0}y^i)+\Psi (-2\alpha
 Qs_0+r_{00})(b^iy^j-b^jy^i)=\frac{1}{2}(G^i_{kl}y^j-G^j_{kl}y^i)y^ky^l,\label{y40}\\
 &&\ \ \ \ \ \ \ \ \    \big(Q:=\frac{\phi'}{\phi-s\phi'},\ \
  \Psi:=\frac{Q'}{2\Delta},\ \
  \Delta:=1+sQ+(b^2-s^2)Q'\big),\nonumber
 \eeq
where $G^i_{kl}$ are some scalar functions. Then plugging
$\phi(s)=se^{\pm 1/s^2}$ into (\ref{y40}), we can prove that
(\ref{y40}) is equivalent to (\ref{y39}). The details are omitted.
\qed

\begin{lem}
 Suppose that (\ref{y4}) holds. Then we have
  \be\label{y040}
 V^c(b^2)=\pm(\tau-2c)b^4.
  \ee
\end{lem}

{\it Proof :}  Rewrite (\ref{y4}) as
 \be\label{y041}
 V_{i|j}+V_{j|i}=2\tau a_{ij}\pm(2c-\tau)b_ib_j,\ \ \ \ V^jb_{i|j}+b^jV_{j|i}=\tau
 b_i.
 \ee
Then we have
 \beqn
 V^c(b^2)&=&2V^jb^ib_{i|j}\\
 &=&2(\tau b^2-b^ib^jV_{i|j}),\ \ \text{(by the second equation of
 (\ref{y041}))}\\
 &=&\pm(\tau-2c)b^4,\ \ \ \ \ \ \ \ \text{(by the first equation of
 (\ref{y041}))}.
 \eeqn
 This completes  the proof.   \qed

 \

Now define a general pair $(\widetilde{\alpha},\widetilde{\beta})$
by
 \be\label{y41}
 \widetilde{\alpha}=\sqrt{u(b^2)\alpha^2+v(b^2)\beta^2},\ \ \widetilde{\beta}=w(b^2)\beta,
 \ee
where  $u=u(t)\ne 0,v=v(t),w=w(t)\ne 0$ are some functions
satisfying
 \be\label{y42}
 u'=uw^{-1}w'\mp t^{-2}u,\ \ \  v'=vw^{-1}w'-t^{-2}(\pm v-u).
 \ee

 \begin{lem}\label{lem59}
Under (\ref{y41}) and (\ref{y42}), the first equation of
(\ref{y39}) implies that $\widetilde{\beta}$ is a conformal 1-form
of $\widetilde{\alpha}$.
 \end{lem}

{\it Proof :} Under (\ref{y41}) and (\ref{y42}), it can be
directly verified
 that the first equation of (\ref{y39}) is reduced to
  $$
 \widetilde{r}_{00}=\pm\frac{\sigma b^2(b^4w'\pm
 w)}{2(u+vb^2)}\ \widetilde{\alpha}^2,
  $$
which completes the proof.  \qed

\begin{lem}
 Under (\ref{y41}), the conditions (\ref{y4}) and (\ref{y42}) give
  \be\label{y44}
 V^c(\widetilde{\alpha}^2)=\Big[(2c+\tau)u\pm\frac{(\tau-2c)b^4uw'}{w}\Big]\widetilde{\alpha}^2,\
 \ \ \ \  V^c(\widetilde{\beta})=\big[ w\tau\pm(\tau-2c)b^4w' \big].
  \ee
\end{lem}

{\it Proof :} Under (\ref{y41}), we first have
 \beqn
 V^c(\widetilde{\alpha}^2)&=&V^c(b^2)\big(u'\alpha^2+v'\beta^2\big)+uV^c(\alpha^2)+vV^c(\beta^2),\\
 V^c(\widetilde{\beta})&=&V^c(b^2)w'\beta+wV^c(\beta).
 \eeqn
Then plugging (\ref{y4}), (\ref{y040}) and (\ref{y42}) into the
above two equations, we obtain (\ref{y44}).   \qed

\begin{lem}
The equation (\ref{y4}) is equivalent to
   \be\label{y37}
 V^c(\widetilde{\alpha}^2)=2\tau \widetilde{\alpha}^2,\ \ \ \ \ \  V^c(\widetilde{\beta})=2(\tau-c)\widetilde{\beta},
  \ee
 where $\widetilde{\alpha}$ and $\widetilde{\beta}$ are given by
 \be\label{y36}
 \widetilde{\alpha}:=\sqrt{\alpha^2+(1-b^{-2})\beta^2},\ \ \ \  \widetilde{\beta}:=e^{\mp
 b^{-2}}\beta.
 \ee
\end{lem}

{\it Proof :} Let
 \be\label{y46}
 u(b^2)=1,\ \ \  v(b^2)=1-\frac{1}{b^2},\ \ \ \  w(b^2)=e^{\mp
 b^{-2}}.
 \ee
It is easy to verify that (\ref{y46}) is a special solution of the
ODE (\ref{y42}). Then plugging (\ref{y46}) into (\ref{y41}) we get
(\ref{y36}). By (\ref{y44}) and (\ref{y46}), we get (\ref{y37}).
\qed

\

 Since $F$ is a Landsberg metric in the dimension $n\ge 3$, or a
Douglas metric, we have (\ref{y39}) by Lemma \ref{lem57}. Then by
Lemma \ref{lem59} we see that $\widetilde{\beta}$ is a conformal
1-form of $\widetilde{\alpha}$. Finally, by Lemma \ref{lem51} and
(\ref{y37}) we see that $c$ is a constant. So $V$ is a homothetic
vector field. This completes the proof of Theorem \ref{th53}. \qed

\vspace{0.6cm}

\noindent Guojun Yang \\
Department of Mathematics \\
Sichuan University \\
Chengdu 610064, P. R. China \\
yangguojun@scu.edu.cn

\end{document}